\documentclass[global,
]{svjour}

\usepackage[latin1]{inputenc}
\usepackage{latexsym}
\usepackage{graphics}
\usepackage{times}
\usepackage{mathrsfs,amsmath,amsbsy,amssymb,amsfonts,txfonts,dsfont}

\newcommand{\pg}[1]{\left\{#1\right\}}
\newcommand{\pq}[1]{\left[#1\right]}
\newcommand{\pt}[1]{\left(#1\right)}
\newcommand{\vir}[1]{``#1''}
\renewcommand{\d}{\mathrm{d}}
\def\1{\mathbf{1}}

\begin{document}
\title{
On asymptotically efficient maximum likelihood
estimation of linear functionals in Laplace measurement error models
}


\author{}


\author{Catia Scricciolo
\thanks{Catia Scricciolo\\ 
Dipartimento di Scienze Economiche, Universit\`a degli Studi di Verona,
Polo Universitario Santa Marta, Via Cantarane 24,
I-37129 Verona (VR), ITALY, \email{catia.scricciolo@univr.it}}}

\institute{}

%

%
\date{
}
%

\titlerunning{Asymptotically efficient MLE}

\authorrunning{C. Scricciolo}

\maketitle

\begin{abstract}
Maximum likelihood estimation of linear functionals in the inverse problem of deconvolution is considered.
Given observations of a random sample from a distribution $P_0\equiv P_{F_0}$ 
indexed by a (potentially infinite-dimensional) parameter $F_0$, which is the distribution of the latent variable in a standard additive Laplace measurement error model, one wants to estimate a linear functional of $F_0$.
Asymptotically efficient maximum likelihood estimation (MLE) of integral linear functionals of the mixing distribution $F_0$ in a convolution model with the Laplace kernel density is investigated. 
Situations are distinguished in which the functional of interest can be consistently estimated at $n^{-1/2}$-rate by the plug-in MLE, 
which is asymptotically normal
and efficient, in the sense of achieving
the
variance lower bound, from those in which 
no integral linear functional can be estimated at parametric rate, which precludes any possibility for asymptotic efficiency. 
The $\sqrt{n}$-convergence 
of the MLE, valid in the case of a degenerate mixing distribution at a single location point, fails 
in general, as does
asymptotic normality. It is shown that
there exists no regular estimator sequence for integral linear functionals of the mixing distribution 
that, when recentered about the estimand and $\sqrt{n}$-rescaled, is asymptotically efficient, 
\emph{viz}., has Gaussian limit distribution with minimum variance. One can thus only expect estimation with some slower rate and,
often, with a non-Gaussian limit distribution.
\end{abstract}

\section{Introduction}
\label{intro}
\allowdisplaybreaks
The problem of asymptotically efficient estimation of integral linear functionals of the distribution of the latent variable in a standard additive Laplace measurement error model is considered. The focus is on establishing whether asymptotic normality and efficiency hold for the estimator obtained by plugging into the functional of interest the NPMLE of the mixing distribution in a convolution model with the Laplace kernel density. 
We study the behaviour of the plug-in NPMLE to answer the question of whether there exist integral linear functionals
of the mixing distribution that can be consistently estimated by the maximum likelihood method at $n^{-1/2}$-rate, 
the recentered and $\sqrt{n}$-rescaled version of the plug-in NPMLE 
being asymptotically normal with zero mean and minimum variance. Situations are distinguished in which the plug-in NPMLE is 
consistent at parametric rate and asymptotically
efficient,
albeit the mixing distribution itself can typically be estimated only at slower rates, from those in which there exists no regular sequence of estimators that can be asymptotically 
efficient. The model is described hereafter and the problem formally stated.


\paragraph{Model description}
Let $X$ be a real-valued random variable (r.v.) with distribution $P_0$
defined, for every Borel set $B$ on the real line,  
by the mapping $B\mapsto P_0(B):=\mathrm{P}(X\in B)$.
Suppose that $P_0$ is dominated by Lebesgue measure $\lambda$ on $\mathbb R$, 
with probability density function (p.d.f.) $p_0:=\d P_0/\d \lambda$. 
Let $X$ satisfy the relationship
\begin{equation}\label{eq:conv}
X=Y+Z,
\end{equation}
where $Y$ and $Z$ are (stochastically) independent, \emph{unobservable} random variables such that $Y$ has \emph{unknown} cumulative distribution function (c.d.f.) $F_0$ and
$Z$ has the \emph{standard classical} Laplace\footnote{It is also known as the \emph{first law of Laplace}
to distinguish it from the \emph{second law of Laplace}, as the normal distribution is sometimes called. It was named after Pierre-Simon Laplace (1749--1827) who, in 1774 (Laplace 1774), obtained 
$e^{-|z-\theta|}/2$, for $z,\,\theta\in\mathbb{R}$,
as the density of the distribution whose likelihood is maximized when the location parameter $\theta$ is equal to the sample median.}
or double exponential\footnote{It is so called because
it is formed by reflecting the exponential distribution around its mean.}
distribution with scale parameter $s=1$, in symbols, $Z\sim\mathrm{Laplace}\,(0,\,1)$, whose density $k$ has expression 
\begin{equation}\label{eq:laplace}
k(z)=\frac{1}{2}e^{-|z|}, \quad z\in\mathbb{R}.
\end{equation}
The density $p_0$ is therefore the convolution of $F_0$ and $k$ or a location mixture of Laplace densities with mixing distribution $F_0$ supported on a subset $\mathscr{Y}\subseteq\mathbb{R}$, 
where $\mathscr Y$ stands for \emph{a} support of $F_0$, see, \emph{e.g.}, Billingsley (1995), p. 23, 
\[p_0(x)\equiv p_{F_0}(x)=\int_{\mathscr{Y}}k(x-y)\,\d F_0(y)=\int_{\mathscr{Y}}\frac{1}{2}e^{-|x-y|}\,\d F_0(y),\quad x\in\mathbb{R}.\]
For ease of exposition, the density of the standard Laplace distribution 
is considered as a kernel, but the density of \emph{any} 
Laplace distribution centered at zero, with \emph{known} scale parameter $s>0$, in symbols, $Z\sim\mathrm{Laplace}\,(0,\,s)$, 
whose variance is $\sigma^2_Z=2s^2$,\footnote{To see that $\sigma^2_Z=2s^2$, one can take into account
that, if $V_1$ and $V_2$ are independent r.v.'s, identically distributed as an exponential with parameter $1/s$, in symbols, $V_j\sim \textrm{Exp}\,(1/s)$, $j=1,\,2$, then $V_1-V_2$ has a $\textrm{Laplace}\,(0,\,s)$ distribution. Consequently, 
$\sigma^2_Z=2\sigma^2_{V_1}=2(1/s)^{-2}=2s^2$.}
could be employed. 
Assume that $X_1,\,\ldots,\,X_n$ constitute a random sample from $p_0$. Every 
$X_i$ then satisfies
\begin{equation}\label{eq:conveq}
X_i=Y_i+Z_i,\quad i=1,\,\ldots,\,n,
\end{equation}
where $Y_1,\,\ldots,\,Y_n$ and $Z_1,\,\ldots,\,Z_n$ are independent samples from the distributions with c.d.f. $F_0$ and p.d.f. $k$, respectively. The r.v.'s $Y_1,\,\ldots,\,Y_n$ are independent and identically distributed (i.i.d.) as $Y$ and $Z_1,\,\ldots,\,Z_n$ are independent copies of $Z$. Realizations 
 $x_1,\,\ldots,\,x_n$ of the noisy sample data $X_1,\,\ldots,\,X_n$ are observed instead of outcomes of the uncorrupted r.v.'s $Y_1,\,\ldots,\,Y_n$. The r.v.'s $Z_1,\,\ldots,\,Z_n$ represent additive errors and their distribution is called the \emph{error distribution}. In this model, the variable of interest $Y$ cannot be directly observed and empirical access is limited to the sum of $Y$ and the \vir{noise} $Z$. Therefore, estimating the distribution function $F_0$ of $Y$, or related quantities like the p.d.f. $f_0$ (if it exists), based on a sample from $P_0$, accounts for solving a particular \emph{inverse} problem, called \emph{deconvolution}, which consists in reconstructing (estimating) $F_0$ from indirect noisy observations $X_1,\,\ldots,\,X_n$ drawn from $P_0\equiv P_{F_0}$, the latter being the image of $F_0$ under a known transformation that has to be \vir{inverted}. 
As remarked in Groeneboom and Wellner (1992),
p. 4, as well as in Bolthausen \emph{et al.} (2002),
p. 363, the problem can be viewed as a \emph{missing data problem}: the
complete observations would consist of the independent pairs $X^0_i:=(Y_i,\,Z_i)$, with $X^0_i\sim Q:=F_0\times F_Z$,\footnote{The symbol $F_Z$ denotes the c.d.f. of $Z$, that is, $F_Z(z)=\int_{-\infty}^z k(u)\,\d u$, with density $k$ as in \eqref{eq:laplace}.} but part of the data is missing and only outcomes or realizations of the sums $X_i=T(X_i^0):=(Y_i+Z_i)\sim P_0\equiv QT^{-1}$, viewed as transformations of the $X^0_i$'s through the function $T$, are observed.


The statistical model described by relationship \eqref{eq:conv}, with a zero-mean Laplace measurement error r.v. $Z$ independent of $Y$, is a special case of the \emph{classical error model} $X=Y+Z$, in which $X$ is a measurement of $Y$ in the usual sense, 
$Z$ has zero mean and is independent of $Y$,
see, \emph{e.g.}, Buzas \emph{et al.} (2005), p. 733, and the references therein.
Measurement errors with possibly different structures occur in nearly every discipline from medical statistics to astronomy and econometrics, 
\emph{cf.} the monographs of 
Fuller (1987)
and Buonaccorsi (2010).
Furthermore, 
the Laplace distribution finds applications in a variety of disciplines, 
from image and speech recognition to ocean engineering, see Kotz \emph{et al.} (2001), chpts. 7--10, pp. 343--397.
An application to quality control of the classical Laplace measurement error model is outlined hereafter.

\paragraph{Application to steam generator inspection}
An application of the Laplace measurement error model to steam generator inspection and testing is described herein, see Easterling (1980) and 
Sollier (2017) for more details. Steam generators of pressurized water reactors contain 
many tubes through which heated water flows. For a variety of reasons, such as corrosion-induced wastage, the steam generator tube integrity can be degraded, the walls becoming thinned or
cracked. Leaks may occur during normal operating conditions, thus requiring the plant to be shut down. 
In order to develop an inspection plan, a statistical model for tube degradation is considered. Experimental data evidentiate that measurements are affected by heavy-tailed and biased errors that can be represented by a r.v. $E$ following a Laplace distribution with mean $\mu>0$ and scale parameter $s>0$, in symbols, $E\sim\mathrm{Laplace}\,(\mu,\,s)$, with density 
\[f_E(e)=\frac{1}{2s}\exp{(-|e-\mu|/s)}, \quad e\in\mathbb R.\]
Denoted by $D$ the \emph{actual} degradation (extent of thinning)
of a tube, expressed as a percentage of the initial tube wall
thickness, the \emph{measured} degradation $M$ is modeled as 
\begin{equation*}\label{eq:degradation}
M=D+E,
\end{equation*}
where $E$ is supposed to be independent of $D$.
Assuming that the scale parameter $s$ is known and the distribution of $D$ possesses probability 
density function, say $f_D$,   
the interest is, in the first place, in estimating the p.d.f. of $D$, based on i.i.d. observations $M_1,\,\ldots,\,M_n$ drawn from the distribution of $M$. An exponential distribution for $D$, 
with scale parameter $\tau>0$, in symbols, $D\sim\textrm{Exp}\,(1/\tau)$, whose density has expression
\[f_D(d)=\frac{1}{\tau}\exp{(-d/\tau)},\quad d>0,\]
provides an exponential-double exponential model for the actual degradation $M$, which has proved to have an adequate fit on experimental data. Statistical procedures for fitness-for-service assessment are described in Carroll (2017).

\paragraph{Asymptotic efficiency of the NPMLE for linear functionals of the mixing distribution}
For many purposes, interest can lie in only 
few aspects of the distribution of $Y$, key features of which can be represented as linear functionals of $F_0$. 
In what follows, symbols $F_0$ and $F$ will be used to indicate probability measures (p.m.'s) on $(\mathscr{Y},\,\mathscr{B}(\mathscr{Y}))$,
where $\mathscr{B}(\mathscr{Y})$ denotes the Borel $\sigma$-field on $\mathscr{Y}$,
as well as the corresponding cumulative distribution functions, the correct meaning being clear from the context.
Letting
\[\mathscr P:=\{\mbox{all p.m.'s $F$ on $(\mathscr{Y},\,\mathscr{B}(\mathscr{Y}))$}\}\]
be the collection of all probability measures $F$ on $(\mathscr{Y},\,\mathscr{B}(\mathscr{Y}))$,
a functional is a mapping $\psi:\,\mathscr P\rightarrow \mathbb R$
that maps every $F\in\mathscr P$ to a real number $\psi(F)$.
The focus is on estimating \emph{integral} linear functionals
\begin{equation}\label{eq:2}
F\mapsto\psi({F}):=\int_{\mathscr{Y}}a(y)\,\d F(y)
\end{equation}
at the \vir{point} $F_0$, where the function $a\in L^1(F_0)$ is given.
The following examples illustrate choices of $a$ 
for some common statistical functionals.\\[-12pt]
\begin{itemize}
\item {{\emph{Distribution function at a point}}} If, for some fixed $y_1\in\mathbb{R}$, the function 
$a(\cdot)=1_{(-\infty,\,y_1]}(\cdot)$, then 
$\psi(F_0)=\int_{\mathscr{Y}}1_{(-\infty,\,y_1]}(y)\,\d F_0(y)=F_0(y_1)$ 
is the c.d.f. of $Y$ at the point $y_1$.\\[-8pt]
\item {{\emph{Probability of an interval}}}
If, for fixed points $y_1,\,y_2\in\mathbb{R}$, the function 
$a(\cdot)=1_{(y_1,\,y_2]}(\cdot)$, then 
$\psi(F_0)=\int_{\mathscr{Y}}1_{(y_1,\,y_2]}(y)\,\d F_0(y)=F_0(y_2)-F_0(y_1)=\mathrm P(y_1<Y\leq y_2)$ is the probability of the interval 
$(y_1,\,y_2]$.\\[-8pt]
\item {{\emph{Mean}}}
If
$a(\cdot)=\textrm{id}_{\mathscr Y}(\cdot)$ is the identity function on $\mathscr Y$, then $\psi(F_0)=\int_{\mathscr{Y}}y\,\d F_0(y)
=\mathrm{E}Y$ is the expected value of $Y$
which, for \emph{any} kernel density $k$ (not necessarily the Laplace) with zero mean, $\mathrm{E}Z=0$, is equal to $\mathrm{E}X$: in fact, from the relationship $X=Y+Z$
in \eqref{eq:conv}, it follows that
$\mathrm{E}X=\mathrm{E}Y+\mathrm{E}Z=\mathrm{E}Y$ by linearity of the expected value.
\\[-8pt]
\item{\emph{$r$th moment}} 
If, for any positive integer $r$, the function
$y\mapsto a(y)=y^r$, then $\psi(F_0)=\int_{\mathscr{Y}}y^r\,\d F_0(y)
=\mathrm{E}Y^r$ is the $r$th moment of $Y$.\\[-8pt]
\item {{\emph{Moment generating function}}} If, for some fixed point $t\in\mathbb R$ such that $0<|t|<1$, the mapping $y\mapsto a(y)=e^{ty}$, then $\psi(F_0)$ coincides with the moment generating function (m.g.f.) of $F_0$ at $t$, denoted by $M_{F_0}(t)$ or $M_Y(t)$, that is, 
$\psi(F_0)=\int_{\mathscr{Y}}e^{ty}\,\d F_0(y)=M_{F_0}(t)$.
\noindent
Some features of the mixing distribution $F_0$, like the mean or the variance, can be 
expressed in terms of the derivatives of the corresponding m.g.f. $M_{F_0}$ evaluated at zero. Therefore, in principle, results for estimating aspects of $F_0$ can be obtained as by-products of the inference on $M_{F_0}$.
\end{itemize}

\smallskip

A standard and principled method for pointwise estimation of linear functionals
consists in plugging \emph{the}\footnote{Uniqueness of $\hat F_n$ is not guaranteed, it is therefore with an abuse of language that we refer to \emph{the} NPMLE throughout the article.} NPMLE $\hat F_n$ of $F_0$ into $\psi(\cdot)$ to obtain the \emph{plug-in estimator}
$$\psi({\hat F_n}):=\psi(F)\big|_{F=\hat F_n}.$$
A NPMLE $\hat F_n$ of $F_0$ is a measurable function of the observations $X_1,\,\ldots,\,X_n$ taking values in 
$\mathscr{P}$, which
is not necessarily uniquely defined by the relationship
\[\hat F_n\in\mathop{\arg\max}_{\substack{F\in\mathscr P}}
\frac{1}{n}\sum_{i=1}^n\log p_F(X_i),\]
equivalently written as
\[\hat F_n\in
\mathop{\arg\max}_{\substack{F\in \mathscr P}}
\mathbb P_n\log p_F,\]
where
$$p_F(\cdot):=\int_{\mathscr{Y}}k(\cdot-y)\,\d F(y)$$ is the generic location mixture of Laplace densities with mixing distribution $F$ and
$${\mathbb P_n}:=\frac{1}{n}\sum_{i=1}^n\delta_{X_i}$$ 
is the empirical probability measure associated with the random sample $X_1,\,\ldots,\,X_n$, namely, the discrete uniform distribution on the sample values that puts mass $1/n$ on each one of the observations. In the sequel, 
for a measurable function $f:\,\mathscr X\rightarrow\mathbb R$, 
where $\mathscr X\subseteq\mathbb R$ is specified at the different occurrences, the notation $\mathbb P_n f$ is used to abbreviate the empirical average $n^{-1}\sum_{i=1}^nf(X_i)$. Analogously, $P_0f$ is used in lieu of $\int f\,\d P_0$.
Hereafter, unless it is necessary within the context to specify the integral domain, integration is understood to be performed over the entire natural domain of the integrand. Throughout the article, the probability measure $\mathbf{P}$ stands for $P_0^n$, the joint law of the first $n$ coordinate projections of the infinite product probability measure $P_0^{\mathbb{N}}$. 
Sequences of random variables are meant to convergence (in law or in probability) as the sample size $n$ grows indefinitely large
(as $n\rightarrow +\infty$).

\paragraph{Historical and conceptual background, overview of the results}
The \emph{deconvolution} problem has been intensively studied over the last thirty years. There exists a vast literature on \emph{density} deconvolution, which accounts for reconstructing/estimating the density $f_0$ of $Y$ (if it exists) that satisfies the equation
\[p_0(x)=\int_{\mathscr Y} k(x-y)f_0(y)\,\d y,\quad x\in\mathbb R,\]
wherein the kernel density $k$ (not necessarily a Laplace) is assumed to be known, based on outcomes $x_1,\,\ldots,\,x_n$ of i.i.d. r.v.'s $X_1,\,\ldots,\,X_n$ as in \eqref{eq:conveq}. We cite key articles of the early 90's like Carroll and Hall (1988),
Stefanski and Carroll (1990),
Fan (1991)
\footnote{For a recent reference list, see also Davidian \emph{et al.} (2014).}
that have been ground-breaking and have had a great impact
on the area of measurement error, setting the general framework for attacking measurement error/deconvolution problems
and developing an approach based on Fourier inversion techniques to construct a deconvolution kernel density estimator for recovering the density of the latent distribution, 
meanwhile showing how difficult it is to account for measurement errors: in fact, the smoother the error distribution, the stronger its confounding effect on the latent distribution, hence, the slower the optimal attainable rate of convergence for its estimators.

Far less instead seems to be known about \emph{distribution function} deconvolution, keynote 
contributions, also based on Fourier inversion techniques, being those of Hall and Lahiri (2008), Dattner \emph{et al.} (2011), the former article containing an illuminating critical analysis of the background to the problem of distribution estimation 
in deconvolution problems. 
Since the focus of this article is on the behaviour of the NPMLE $\hat F_n$ of 
the mixing distribution 
$F_0$, attention is hereafter restricted to review the theory of non-parametric maximum likelihood estimation in deconvolution problems. 
In general mixture models, the NPMLE $\hat F_n$ of $F_0$ is \emph{discrete}, with at most $l\leq n$ support points, $l$ being the number of distinct 
values of the data points, \emph{cf.} Lindsay (1983).
In deconvolution problems with continuous and symmetric (about the origin) kernels decreasing on $[0,\,+\infty)$, \emph{a} NPMLE $\hat F_n$ always exist (uniqueness is not guaranteed), see Groeneboom and Wellner (1992), 
Lemma 2.1, pp. 57--58; for kernels that are also strictly convex on $[0,\,+\infty)$, like the (standard) Laplace, 
the NPMLE $\hat F_n$ is supported on the set of observation points $\{X_1,\,\ldots,\,X_n\}$,
so that the corresponding probability measure, still denoted by $\hat F_n$ consistently with the notational convention adopted throughout, is concentrated on the range of the data points $[X_{1:n},\,X_{n:n}]$,\footnote{Following a common notational convention, we denote by
$X_{r:n}$ the $r$th order statistic.}
$$\hat F_n([X_{1:n},\,X_{n:n}])=1,$$
see \emph{ibid.},
Corollary 2.2 and Corollary 2.3, p. 59 and p. 60, respectively. Consistency of $\hat F_n$ at a \emph{continuous} distribution function $F_0$ is proved in \emph{ibid.}, 
\S\,4.2, pp. 79--81. More is known about the \emph{one}-parameter (location only) Laplace model, which can be viewed as a degenerate mixture with point mass mixing distribution at some fixed $\theta\in\mathbb R$, the Dirac measure at $\theta$.\footnote{The Dirac measure at $\theta$, denoted by $\delta_\theta(\cdot)$, is defined on the Borel sets $B\in\mathscr B(\mathbb R)$ by $\delta_\theta(B)=1,\,0$ if $B\ni\theta$ or $B\notni \theta$, respectively.} A simple maximization argument to find a MLE of the location parameter $\theta$, denoted by $\hat\theta_n$, is given in Norton (1984): the \emph{sample median} is a MLE which is an \emph{M-estimator}, see Huber (1967), solving the equation $\sum_{i=1}^n\mathrm{sign}(X_i-\theta)=0$.\footnote{The \emph{sign-function} is defined as $\mathrm{sign}(x)=-1,\,0,\,1$ if $x<0$, $x=0$ or $x>0$, respectively.}
Therefore, 
a MLE exists, but may not be unique: if $n$ is odd, that is, $n=2m+1$ for some $m\in\mathbb N$, then the sample median is uniquely defined as the middle observation $X_{m+1:n}$, while, if $n=2m$ is even, then 
there are two middle observations $X_{m:n}$ and $X_{m+1:n}$ 
so that, in principle, any value in the interval $[X_{m:n},\, X_{m+1:n}]$ 
could be chosen, even if the \emph{canonical median}
$(X_{m:n}+X_{m+1:n})/2$, the average of the middle observations, is typically used in practice. Therefore,
\[\hat\theta_n=
  \begin{cases}
    X_{m+1:n}, &\,\,\,\text{for } n=2m+1,\\[-0.25cm]
 &\qquad\qquad\qquad\qquad\text{} m\in\mathbb N,\\[-0.25cm]
    \frac{1}{2}(X_{m:n}+X_{m+1:n}),  &\,\,\,\text{for } n=2m.
  \end{cases}
\]
The sample median $\hat\theta_n$ is a MLE and is \emph{asymptotically efficient}, that is, 
consistent and, when recentered at 
$\theta$ and $\sqrt{n}$-rescaled, asymptotically normal with zero mean and variance equal to one, which is the information lower bound corresponding to the amount of information in a single observation,
\[\sqrt{n}\big(\hat\theta_n-\theta\big)\xrightarrow{\mathscr{L}}\,\mathscr{N}\big(0,\,1\big),\]
where \vir{$\,\xrightarrow{\mathscr{L}}\,$} denotes convergence in law.
This has been established by Daniels (1961) 
who, motivated by the non-differentiability at zero of the (standard) Laplace density, proved a general theorem on the asymptotic efficiency of the MLE under conditions not involving the second and higher-order derivatives of the likelihood function. Even though, as later on noted by Huber (1967), a crucial step has been overlooked in Daniels' proof, 
the assertion remains valid.

In this article, we study the behaviour of the plug-in NPMLE $\psi(\hat F_n)$ to answer the question of whether there exist integral linear functionals $\psi(F_0)$ that can be consistently estimated by the maximum likelihood method at $n^{-1/2}$-rate and for which $$\sqrt{n}\big(\psi(\hat F_n)-\psi(F_0)\big)$$ is \emph{asymptotically efficient}, in the sense of Definition 2.8 in 
Bolthausen \emph{et al}. (2002), p. 349, that is, asymptotically normal with zero mean and minimum variance. 
In fact, also in non-parametric problems estimation can be performed at $n^{-1/2}$-rate and, in general mixture models 
$\int_{\mathscr Y} k(\cdot\mid y)\,\d F_0(y)$, where $k$ is any (not necessarily the Laplace) kernel density, there may exist linear functionals of $F_0$ that are estimable at parametric rate, even if $F_0$ itself can be pointwise estimated only at slower rates. 
Fundamental contributions developing the theory of information bounds are van der Vaart (1991), van der Vaart (1998), chpt. 25, pp. 358--432, Bolthausen \emph{et al.} (2002), Part III, pp. 331--457,
Groeneboom and Wellner (1992), with emphasis on non-parametric maximum likelihood estimation, and
van de Geer (2000), chpt. 11, pp. 211--246, with a focus on asymptotic efficiency of the NPMLE in mixture models.
To exemplify the issue, consider estimating the mean functional $\psi(F_0)=\mathrm{E}Y$ which, in a mixture model 
such that $\mathrm E(X\mid Y)=Y$, is equal to $\mathrm E X$. 
Then, the sample mean $\bar X_n=n^{-1}\sum_{i=1}^n X_i$ is a $n^{-1/2}$-consistent and, after $\sqrt{n}$-rescaling, asymptotically normal estimator of $\mathrm{E}Y$, 
but may not be a MLE; furthermore, it does not take into account the information that the sampling density is a mixture. 
On other side, the MLE may be $n^{-1/2}$-consistent and converge to a normal distribution with smaller variance than that of the sample mean. 
This is the case for the sample median in the single-parameter (location only) Laplace model, see van der Vaart (1998), Example 7.8 on location models, p. 96. 
Surprisingly, little is known in general about the asymptotic behaviour of the plug-in NPMLE for linear functionals in Laplace convolution mixtures,
even if only for estimating the mean functional. Although the topic is useful, existence of this gap can be partially explained by the fact that the 
Laplace or double-exponential distribution is not an exponential family model so that standard results may not be valid or immediately available from the theory of exponential families.


In order to investigate whether integral linear functionals of the mixing distribution in a convolution 
model with the Laplace kernel density are estimable at $n^{-1/2}$-rate, we appeal to van der Vaart's differentiability theorem, 
see Theorem 3.1 in van der Vaart (1991), p. 183, a general result that allows for a unified treatment of the information lower bound theory 
based on the concept of a differentiable functional, see, for a definition of the latter, display (2.2) in van der Vaart (1991), p. 180, or Definition 1.10 in 
Bolthausen \emph{et al.} (2002), p. 343.
The differentiability theorem characterizes differentiable functionals and, 
by combining the description of the set of differentiable functionals with a result stating that the existence of \emph{regular} estimator sequences for a functional implies its differentiability, provides a way to distinguish situations in which the functional of interest is estimable at $n^{-1/2}$-rate from situations in which this is not the case, see van der Vaart (1998), p. 365, for a definition of a regular estimator sequence.
A necessary and sufficient condition for differentiability of a (not necessarily linear) functional is that its gradients are contained in the range of the adjoint of the score operator, where the score operator can be viewed as a derivative (in quadratic mean) of the map $F\mapsto P_F$,
see (3.6) in van der Vaart (1991), p. 183, or (25.29) in \S\,25.5 of van der Vaart (1998), p. 372.
As previously mentioned, differentiability is necessary for regular estimability of a functional or, equivalently, for
the existence of regular estimator sequences, see Theorem 2.1 of van der Vaart (1991), p. 181,
so that if the functional is not differentiable, then there exists no regular estimators
and estimation at $n^{-1/2}$-rate is impossible.
Interestingly, for real-valued functionals, the differentiability condition is equivalent to having positive efficient information, see Theorem 4.1 of van der Vaart (1991), pp. 186--187. We find that the differentiability condition fails for integral linear functionals of the mixing distribution in a convolution model with the Laplace density, this implying that there exists no estimator sequence for $\psi(F_0)$ that is regular at $F_0$ and estimation at $n^{-1/2}$-rate is impossible.

\paragraph{Organization}
The rest of the article is organized as follows. 
The main results are presented in Sect. \ref{section:1}, 
which is split into two parts. In the first one, asymptotic efficiency of the plug-in NPMLE for integral linear functionals of the mixing distribution in a convolution model with the (one-sided) exponential kernel density 
is 
analysed and set in the affirmative. Construction of interval estimators and tests based on a Studentized version of the plug-in NPMLE, when the asymptotic variance is consistently estimated, is revisited. Conditions for extending results to non-linear functionals are discussed as a side-issue. In the second part, the focus is on asymptotically efficient estimation by the plug-in NPMLE for integral linear functionals of the mixing distribution in a convolution model with the 
double-exponential (Laplace) kernel density. It is shown that, except for the case of a degenerate mixing distribution at a single location point, maximum likelihood estimation completely fails, in the sense that no integral linear functional can be estimated at $n^{-1/2}$-rate, which precludes any possibility for the NPMLE of being asymptotically efficient. Indeed, 
there exists no regular sequence of estimators for integral linear functionals of the mixing distribution
that can be asymptotically efficient, therefore, estimation of linear functionals  is impossible at $n^{-1/2}$-rate.
Final remarks and comments are exposed in Sect. \ref{sec:3}. 
Proofs of the main results are deferred to the appendices: 
Appendix A reports the proof of the result for convolution with the exponential density, Appendix B
reports the proof of the result for convolution with the Laplace density.




\section{Main results}\label{section:1}
In this section, the main results of the article are presented.
First, the case of a convolution model with an exponential kernel density is considered. 
Since, as previously noted, the Laplace density in \eqref{eq:laplace} can be thought of as two exponential densities spliced together back-to-back, the positive half being a
standard exponential density scaled by $1/2$, it is reasonable to begin the analysis from the problem of asymptotically efficient non-parametric maximum likelihood estimation of linear functionals in a convolution model with the \emph{exponential} kernel density. A preliminary study of the one-sided problem, beyond being of interest in itself, is useful to attack the   
two-sided one by partially reducing it to the previously solved case; it may furthermore provide insight for a better understanding of the reasons why symmetrization leads to a failure of asymptotically efficient estimation of linear functionals in the \emph{double-exponential} (Laplace) case.

\paragraph{Convolution with the exponential density}
In this paragraph, 
a standard exponential kernel density on $[0,\,+\infty)$ is considered. This gives rise to a one-sided mixture density generating the data,
\[p_0(x)\equiv p_{F_0}(x)=\int_{\mathscr Y} e^{-(x-y)}1\{y\leq x\}\,\d F_0(y),\quad x\in\mathscr X,\]
where
$\mathscr Y:=\mathrm{support}\,(F_0)
$
is assumed, without loss of generality, to be
a proper, left-closed 
subset of the 
real line, 
and $\mathscr X:=[y_{\textrm{min}},\,+\infty)$ is the support of $P_0$, 
with $y_{\textrm{min}}:=\mathrm{min}\,\mathscr Y>-\infty$. 
Proposition \ref{prop:1} below
establishes that, under sufficient conditions,
an integral linear functional $\psi(F_0)$
can be consistently estimated at $n^{-1/2}$-rate by the 
plug-in NPMLE $\psi(\hat F_n)$,
which, when recentered about the estimand and $\sqrt{n}$-rescaled, is asymptotically normal and efficient. 
In stating hereafter the assumptions, Newton's notation (or the dot notation) for differentiation is adopted, 
that is, ${\dot a}(y):=\d a(y)/\d y$.

\medskip

\noindent
{{\bf{{Assumptions}}}}
\begin{enumerate}
\item[$(\bf{A0})$] $a\in L^1(F_0)$,\\[-8pt]
\item[$(\bf{A1})$] $M_{F_0}(1):=\int_{\mathscr Y}e^{y}\,\d F_0(y)<+\infty$,\\[-8pt]
\item[$(\bf{A2})$] $|\psi(\hat F_n)-\psi(F_0)|=o_{\mathbf{P}}(1)$,\\[-8pt] 
\item[\quad$(\bf{A3})$] 
(i) \,\,\,\,$a$ is continuous on $\mathscr Y$,\\[1pt]
(ii)\,\, either $\mathscr Y$ is compact or $a$ is bounded on $\mathscr Y$,\\[1pt](iii) $F_0(y_{\textrm{min}})=0$,\\[-8pt]
\item[$(\bf{A4})$] there exists ${\dot a}$ on $\mathscr Y$
and 
$
\sup_{y\in\mathscr Y}|\dot a(y)|<+\infty$,
\item[$(\bf{A5})$] there exists a constant $0< c_0<+\infty$ such that 
$$\sup_{y\in\mathscr Y}\bigg|\frac{\d\big(\dot a(y)e^{-y}\big)}{\d F_0(y)}\bigg|\leq c_0.
$$
\end{enumerate}
\smallskip
\noindent
Some remarks and comments on the above listed assumptions are in order. Except for Assumption $({\bf{A2}})$, which concerns the plug-in NPMLE $\psi(\hat F_n)$, all other assumptions involve the function $a$ and/or the mixing distribution $F_0$ that jointly define the functional $\psi(F_0)$.
Specifically, Assumption $(\bf{A0})$ guarantees that $\psi(F_0)$ is well defined. Assumption $(\bf{A1})$ ensures the existence of the moment generating function of $F_0$ at the point $t=1$.
Assumption $(\bf{A2})$ requires consistency of $\psi(\hat F_n)$ at $\psi(F_0)$. 
If $\hat F_n$ converges weakly to $F_0$ in $P_0^n$-probability, then parts (i) and (ii) of
Assumption $(\bf{A3})$ together imply Assumption $(\bf{A2})$
because $a$ is continuous and bounded on $\mathscr Y$.
Sufficient conditions for $\hat F_n$ to converge weakly to $F_0$ in the convolution model with the standard exponential kernel density on $[0,\,+\infty)$ are stated in Groeneboom and Wellner (1992), p. 86; see also Theorem 2.3 of Chen (2017), 
p. 54, for sufficient conditions for strong consistency of
$\hat F_n$ in general mixture models. 
Part (ii) of Assumption $(\bf{A3})$ postulates that either $\mathscr Y$ is a closed and bounded interval
$[y_{\mathrm{min}},\,y_{\mathrm{max}}
]$ or $\mathscr Y$ is a right-unbounded interval $[y_{\mathrm{min}},\,+\infty)$ and $a$ is bounded.
Assumption $(\bf{A4})$ requires $a$ to be differentiable and bounded on $\mathscr Y$, which, in particular, accounts for 
$a$ to be right-differentiable at $y_{\mathrm{min}}$, that is,
$\dot a_+(y_{\mathrm{min}})<+\infty$,\footnote{The right-derivative of $a$ at $y_{\mathrm{min}}$, denoted by $\dot a_+(y_{\mathrm{min}})$, is defined as the one-sided limit   
$\lim_{y\rightarrow y_{\mathrm{min}}+}
[a(y)-a(y_{\mathrm{min}})]/(y-y_{\mathrm{min}})$ if it exists as a real number.} 
and, in the case where $\mathscr Y=[y_{\mathrm{min}},\,y_{\mathrm{max}}]$,
to be also left-differentiable at $y_{\mathrm{max}}$, that is,  
$\dot a_-(y_{\mathrm{max}})<+\infty$.\footnote{The left-derivative of 
$a$ at $y_{\mathrm{max}}$, denoted by $\dot a_-(y_{\mathrm{max}})$, is defined as the one-sided limit $\lim_{y\rightarrow y_{\mathrm{max}}-
}[a(y_{\mathrm{max}})-a(y)]/(y_{\mathrm{max}}-y)$ if it exists as a real number.}
Assumption $(\bf{A5})$ plays its role in the proof of Proposition \ref{prop:1} 
when bounding the worst possible sub-directions, see 
\eqref{eq:suph} in the Appendix A.

\begin{proposition}\label{prop:1}
Under Assumptions \emph{(}{\bf{A0}}\emph{)}--\emph{(}{\bf{A5}}\emph{)}, we have
\begin{equation}\label{eq:conv1}
\sqrt{n}\big(\psi(\hat F_n)-\psi(F_0)\big)=\sqrt{n}
\mathbb P_n b_{F_0}+o_{\mathbf P}(1)
\,\,\bigg(\xrightarrow{\mathscr{L}}\,\mathscr{N}\big(0,\,\|b_{F_0}\|_{2,P_0}^2\big)\bigg),
\end{equation}
where the mapping
$b_{F_0}:\,\mathscr X\rightarrow \mathbb R$, defined as
\begin{equation}\label{eq:effinfcurv}
x\mapsto b_{F_0}(x):=a(x)-\dot a(x)-\psi(F_0),
\end{equation}
is the efficient influence function, whose
squared $L^2(P_0)$-norm 
\begin{equation}\label{eq:var}
\|b_{F_0}\|_{2,P_0}^2:=\int b_{F_0}^2\,\d P_0
\end{equation}
is the efficient asymptotic variance.
\end{proposition}

\smallskip

Proposition \ref{prop:1} establishes that, under sufficient conditions listed as
Assumptions ({\bf{A0}})--({\bf{A5}}), \emph{cf.} Lemma 4.6 of van de Geer (2003), p. 461, and van der Geer (2000), p. 231, 
the plug-in NPMLE $\psi(\hat F_n)$ consistently estimates $\psi(F_0)$ at $n^{-1/2}$-rate; 
furthermore, when recentered at $\psi(F_0)$ 
and $\sqrt{n}$-rescaled, it is asymptotically distributed as a zero-mean Gaussian, with variance
attaining the lower bound given by the squared $L^2(P_0)$-norm of the efficient influence function, which
plays here the same role as the normalized score function for the case of independent sampling from a parametric model 
$\{p_\theta$, $\theta\in\Theta\subseteq \mathbb{R}^d\}$, $d\geq 1$,
\begin{equation*}\label{eq:norminf}
I^{-1}_\theta\dot{\ell}_\theta,
\end{equation*}
where
$\dot{\ell}_\theta(\cdot)=\partial[\log p_\theta(\cdot)]/\partial \theta$ is the \emph{score function} of the model and $I_\theta=\mathrm{E}[\dot{\ell}_\theta\dot{\ell}_\theta^{\mathrm T}]$ the Fisher information matrix for $\theta$.
In a parametric set-up, the minimum variance lower bound reduces to the
Cram\'er-Rao bound, which states that
the inverse of the Fisher information matrix  $I^{-1}_\theta$ is a lower bound on the variance of any $\sqrt{n}$-rescaled unbiased estimator $T_n\equiv T_n(X_1,\,\ldots,\,X_n)$ of $\theta$, in symbols, $\textrm{var}(\sqrt{n}T_n)\geq I^{-1}_\theta$. Therefore, the counterpart of 
$\|b_{F_0}\|_{2,P_0}^2=P_0b_{F_0}^2$ is, by symmetry of ($I_\theta$ hence of) $I^{-1}_\theta$,
\[\mathrm E[(I^{-1}_\theta\dot{\ell}_\theta)(I^{-1}_\theta\dot{\ell}_\theta)^{\mathrm T}]=I_\theta^{-1}\mathrm E[\dot \ell_\theta\dot \ell_\theta
^{\mathrm T}]I_\theta^{-1}=I_\theta^{-1}I_\theta I_\theta^{-1}=I_\theta^{-1}.\]
In general, considered a function $\psi$ that maps $\Theta$ into $\mathbb R^m$, $m\geq1$, and denoted by $\dot \psi_\theta$ the derivative of 
$\theta\mapsto \psi(\theta)$, the matrix $\dot \psi_\theta I^{-1}_\theta\dot\psi_\theta^{\mathrm T}$ is a lower bound on the variance of any $\sqrt{n}$-rescaled unbiased estimator 
of $\psi(\theta)$.

Even if a statement of the result in Proposition \ref{prop:1} appears in Lemma 4.6 of van der Geer (2003), p. 461, as far as we are aware, a complete derivation of the assertion is not available in the literature, \emph{cf.} also van der  Geer (2000), p. 231, so the proof reported in the Appendix A might prove helpful. The underlying idea is outlined hereafter.
A NPMLE 
$\hat F_n$ 
solves the likelihood equation
for every 
path $t\mapsto F_t$, with $\d F_t:= (1+th_F)\,\d F$, starting at a fixed point $F$ corresponding to $t=0$ and direction $h_F$ such that $\int h_F\,\d F=0$, that is, for every parametric sub-model which passes (at $t=0$) through it. In symbols,
\begin{equation}\label{eq:score}
\frac{\d}{\d t}\mathbb P_n \log p_{\hat F_{n,t}}\bigg |_{t=0}=0.
\end{equation}
For ease of notation, let
$$A_{F}h_F(x):=\int_{\mathscr Y} h_F(y)\frac{e^{-(x-y)}}{ p_F(x)}1\{y\leq x\}\,\d F(y).$$
Equation \eqref{eq:score} reduces to $\mathbb P_n b_{\hat F_n}=0$, where
$b_{\hat F_n}(\cdot)=A_{\hat F_n}h_{\hat F_n}(\cdot)$
is the score function (at $t=0$) in an \vir{information loss model}, see, \emph{e.g.}, \S\,25.5.2 in van der Vaart (1998), pp. 374--375.
If $\hat F_n$ dominates $F_0$, which, however, is seldom true, then $-P_0b_{\hat F_n}=\psi(\hat F_n)-\psi(F_0)$ so that 
$\psi(\hat F_n)-\psi(F_0)=(\mathbb P_n -P_0)b_{\hat F_n}$. Asymptotic equicontinuity arguments then yield that 
$(\mathbb P_n -P_0)b_{\hat F_n}=(\mathbb P_n -P_0)b_{F_0}+o_{\bf P}(n^{-1/2})=\mathbb P_n b_{F_0}+o_{\bf P}(n^{-1/2})$ because $P_0b_{F_0}=0$, namely, the score has zero mean. 
So, 
$\sqrt{n}(\psi(\hat F_n)-\psi(F_0))=\sqrt{n}\mathbb P_n b_{F_0}+o_{\bf P}(1)$.
Asymptotic normality follows. The reader is referred to Sect. 11.2 of van de Geer (2000), pp. 211--246, for a more comprehensive treatment of the topic taking into account technical difficulties to which it cannot be here dedicated the necessary space.  

\smallskip

\begin{remark}
For simplicity, a convolution model with an exponential kernel density on $[0,\,+\infty)$ having intensity $\lambda=1$ 
has been considered 
(we warn the reader of the clash of notation with the symbol $\lambda$ previously used to denote Lebesgue measure on the real line), 
but, as revealed by an inspection of the proof of Proposition \ref{prop:1}, the assertion holds true for every $\lambda>0$.
\end{remark}

\begin{remark}
Part (i) of Assumption ({\bf{A3}}) requires $a$ to be continuous  
on $\mathscr Y$, which is not true for indicator functions, 
therefore Proposition \ref{prop:1} does not apply to pointwise estimation 
of the c.d.f. $F_0$ nor to the estimation of the probability of an interval,  
so that it cannot be concluded that these functionals 
are estimable at $n^{-1/2}$-rate by the corresponding plug-in NPMLE's. Indeed, $F_0$ can be \emph{pointwise} estimated only at $n^{-1/3}$-rate, see Groeneboom and Wellner (1992), p. 121. 
Part (ii) of Assumption ({\bf{A3}}) and Assumption ({\bf{A4}})
require that both $a$ and $\dot a$ are bounded on $\mathscr Y$, which, for example, may not be true for the functions $y^r$ and $e^{ty}$ that define the $r$th moment and the moment generating function of $Y$ at the point $t$, respectively: in fact, both $y^r$ and $e^{ty}$, as well as their first derivatives, are continuous on the half-line $[y_{\mathrm{min}},\,+\infty)$, but not bounded therein. Nonetheless, boundedness can be retrieved by restriction to a compact domain.
Therefore, if, besides Assumptions ({\bf{A0}})--({\bf{A2}}) and ({\bf{A5}}), it also holds that $F_0$ has compact support, then, by Proposition \ref{prop:1}, it can be concluded that $\mathrm E Y^r$ and $M_{F_0}(t)$ are consistently and efficiently estimated at $n^{-1/2}$-rate 
by their respective plug-in NPMLE's.
\end{remark}

\smallskip

Although Proposition \ref{prop:1} asserts that certain integral linear functionals can be consistently estimated at $n^{-1/2}$-rate by the 
plug-in NPMLE $\psi(\hat F_n)$, 
which, when recentered at $\psi(F_0)$ and $\sqrt{n}$-rescaled, is asymptotically normal and efficient, two orders of problems may arise that can make it difficult to employ the result for statistical inference: 
\begin{description}
\item[a)] computation of the NPMLE $\hat F_n$,
\item[b)] dependence of the variance $\|b_{F_0}\|_{2,P_0}^2$ on the unknown sampling distribution $P_0$. 
\end{description}
As for the former difficulty, although the NPMLE $\hat F_n$ can be found by a one-step procedure computing the slope of the convex minorant of a certain function, \emph{cf.} Groeneboom and Wellner (1992), pp. 62--63 (see also Vardi (1989) for a different approach), as a by-product of Theorem 11.8 of van de Geer (2000), p. 217, which the assertion of Proposition \ref{prop:1} relies on, the recentered and $\sqrt{n}$-rescaled plug-in NPMLE $\sqrt{n}(\psi(\hat F_n)-\psi(F_0))$ is equivalent, in the sense of being asymptotically approximable, up to an $o_{\mathbf{P}}(1)$-error term, by 
the empirical average of the efficient influence function.
This is part of a general issue concerning the fact that sequences of efficient estimators for functionals are asymptotically approximable by an empirical average of the efficient influence function, 
see, \emph{e.g.}, Lemma 2.9 in Bolthausen \emph{et al.} (2002), p. 349.
In fact, set the position 
\begin{equation}\label{eq:estim}
\tilde\psi_n:=
\frac{1}{n}\sum_{i=1}^n[a(X_i)-\dot a(X_i)]=
\mathbb P_n(a-\dot a)
\end{equation}
and noted that, from the definition of $b_{F_0}$ in \eqref{eq:effinfcurv}, 
the term $\mathbb P_nb_{F_0}$ appearing in \eqref{eq:conv1} writes as 
$\tilde \psi_n-\psi(F_0)$, 
we have
\[\begin{split}
\sqrt{n}\big(\psi(\hat F_n)-\psi(F_0)\big)&=\sqrt{n}\mathbb P_nb_{F_0}+o_{\mathbf P}(1)\\
&=\sqrt{n}\big(\tilde \psi_n-\psi(F_0)\big)+o_{\mathbf P}(1)\,\,\bigg(\,\xrightarrow{\mathscr{L}}\,\mathscr{N}\big(0,\,\|b_{F_0}\|_{2,P_0}^2\big)\bigg).
\end{split}
\]
Thus, both $\sqrt{n}(\psi(\hat F_n)-\psi(F_0))$ and $\sqrt{n}(\tilde \psi_n-\psi(F_0))$
are asymptotically normal and efficient. 
Moreover, estimators arising from $\tilde\psi_n$ may coincide with simple \emph{na\"ive} estimators. For example,\\[-13pt]
\begin{itemize}
\item if $\psi(F_0)=\mathrm{E}Y$, from \eqref{eq:estim}, for
$(a-\dot a)(y)=y-1$, we get the estimator $\tilde\psi_n=\bar X_n-1$, which is the one 
we would suggest considering that $\mathrm{E}Y=\mathrm{E}X-\mathrm{E}Z
=\mathrm{E}X-1$;\\[-7pt]
\item if $\psi(F_0)=M_Y(t)$ for any fixed $t<1$ such that 
$\int_{\mathscr Y} e^{ty}\,\d F_0(y)<+\infty$, then the estimator derived from \eqref{eq:estim}, for
$(a-\dot a)(y)=(1-t)e^{ty}$, is 
$\tilde\psi_n=(1-t)\times{n}^{-1}\sum_{i=1}^n e^{tX_i}$, 
which is the one we would suggest taking into account that $M_Y(t)=M_X(t)/M_Z(t)=(1-t)M_X(t)$, 
where 
$M_{Z}(t)=(1-t)^{-1}$, $t<1$, is the m.g.f. of a standard exponential r.v. $Z$. So, letting $M_n(t):=n^{-1}\sum_{i=1}^n e^{tX_i}$, $t\in\mathbb{R}$, be the empirical m.g.f. for the random sample $X_1,\,\ldots,\,X_n$,
it turns out that 
$\tilde\psi_n={M_n(t)}/{M_Z(t)}$.
\end{itemize}

\smallskip

\noindent
As for the difficulty listed in b), 
by the plug-in approach, replacing the asymptotic variance $\|b_{F_0}\|_{2,P_0}^2$ with a consistent estimator $S_n^2$ leads to the following assertion.
\begin{corollary}\label{cor:1}
Under the conditions of Proposition \ref{prop:1}, 
if, in addition, $S_n^2\,\xrightarrow{\mathrm{P}}\,\|b_{F_0}\|_{2,P_0}^2$, 
where \vir{$\xrightarrow{\mathrm{P}}$} denotes convergence in $P_0^n$-probability,
then
$$
\sqrt{n}\frac{\tilde\psi_n-\psi(F_0)}{S_n}\,\xrightarrow{\mathscr{L}}\,\mathscr{N}(0,\,1).$$
\end{corollary}

\smallskip

\noindent
Replaced the efficient asymptotic variance in \eqref{eq:var} with a consistent sequence of estimators, asymptotic normality of the Studentized version of $\tilde\psi_n$ allows to carry out \emph{pointwise} inference on linear functionals by interval estimation or hypotheses testing constructing confidence intervals or tests, respectively. For every $0<\alpha<1$, let $z_{\alpha/2}$ be the $(1-\alpha/2)$-quantile of a standard normal distribution, \emph{i.e.}, $\Phi(z_{\alpha/2})=1-\alpha/2$, where $\Phi(\cdot)$ stands for the c.d.f. of a standard normal. Then, 
\[C_n:=\big[\tilde\psi_n-z_{\alpha/2} S_n/\sqrt{n},\,\,\,\,\,\,\,\tilde\psi_n+z_{\alpha/2} S_n/\sqrt{n}\big],\quad \mbox{with $\,\mathbf{P}\big(C_n\ni \psi(F_0)\big)=
1-\alpha+o(1)$,}\]
is 
an approximate $(1-\alpha)$-level confidence interval for $\psi(F_0)$.


\begin{remark}
Asymptotic normality of the plug-in NPMLE for linear functionals of the mixing distribution 
can be employed to establish asymptotic normality for \emph{non}-linear functionals. 
Suppose, for instance, that 
$F\mapsto\varphi(F)$ is defined as
\[\varphi(F):=(g\circ\psi)(F)=g\big(\psi(F)\big),\]
where the function $g:\,\mathbb R\rightarrow\mathbb R$ 
has 
non-zero derivative at $\psi(F_0)$ denoted by $\dot g(\psi(F_0))$. 
Asymptotic normality of $\sqrt{n}(\varphi (\hat F_n)-\varphi(F_0))$ then follows from
asymptotic normality of $\sqrt{n}(\psi(\hat F_n)-\psi(F_0))$ by the delta method, see, \emph{e.g.}, chpt. 3 in van der Vaart (1998), pp. 25--34.
So, if the convergence in \eqref{eq:conv1} takes place, then
\[\sqrt{n}\big(\varphi (\hat F_n)-\varphi(F_0)\big)\,\xrightarrow{\mathscr{L}}\,\mathscr{N}\big(0,\,\sigma_{\varphi}^2\big),\quad
\mbox{with  $\sigma_\varphi:=\dot g(\psi(F_0))
\|b_{F_0}\|_{2,P_0}$},\]
where efficiency of $\psi(\hat F_n)$ carries over into efficiency of $\varphi (\hat F_n)$, see \emph{ibid.}, p. 386, for details.

\smallskip

\noindent
Alternatively, set the position
\[R_n:=\frac{\varphi (\hat F_n)-\varphi(F_0)}{\psi(\hat F_n)-\psi(F_0)},
\]
under the condition
$$|R_n-1|=o_{\mathbf{P}}(1),$$
which 
requires that, in probability, 
$\varphi(\hat F_n)-\varphi(F_0)$ behaves asymptotically as $\psi(\hat F_n)-\psi(F_0)$, 
 after $\sqrt{n}$-rescaling, the two differences have the same limiting distribution.
In fact, if the convergence in \eqref{eq:conv1} takes place, then Slutsky's lemma implies that  
\[
\begin{split}
\sqrt{n}\big(\varphi(\hat F_n)-\varphi(F_0)\big)&=
[(R_n-1)+1]
\sqrt{n}\big(\psi(\hat F_n)-\psi(F_0)\big)\\
&=(o_{\mathbf{P}}(1)+1)\sqrt{n}\big(\psi(\hat F_n)-\psi(F_0)\big)
\,\xrightarrow{\mathscr{L}}\,\mathscr{N}\big(0,\,
\|b_{F_0}\|_{2,P_0}^2\big),
\end{split}
\]
see also the Remark of van de Geer (2000) on 
p. 223.
\end{remark}


\paragraph{Convolution with the double-exponential (Laplace) density}\label{paragraph:Laplace}
In this paragraph, the case of main interest of the article concerning asymptotically 
efficient maximum likelihood estimation of linear functionals of the mixing distribution
in a convolution model with the (standard) Laplace kernel density is considered. It has been recalled in Sect.~\ref{intro} 
that, for a \emph{one}-parameter $\theta$ (location only) Laplace model, the sample median $\hat\theta_n$ is a MLE, consistent and asymptotically efficient, even if, for small sample sizes, it may not be the best estimator to use because there exist other unbiased estimators with smaller variances, which are therefore more efficient, see, \emph{e.g.}, Remark 2.6.2 in Kotz \emph{et al.} (2001), p. 82. 
More precisely, for a sample of odd size $n$ from a general $\mathrm{Laplace}\,(\theta,\,s)$ distribution, the 
variance of $\hat\theta_n$ is equal to
$$
\frac{1}{4(n+2)[k(0)/s]^2}=
\frac{s^2}{(n+2)}
,$$ 
where $k(\cdot)$ is the density of a standard Laplace distribution as defined in \eqref{eq:laplace},
while the asymptotic variance is equal to
$$\frac{1}{4n[k(0)/s]^2}=
\frac{s^2}{n}.$$ 
It is just the case to observe that also the sample mean $\bar X_n$ 
is asymptotically normal with mean $\theta$, but the asymptotic relative efficiency (ARE) of the median to the mean, namely, the
ratio of the variance of the sample mean
to the asymptotic variance of the sample median equals $2$:
$$\frac{2s^2/n}{s^2/n}=2.
$$
On a side note, we recall that, for any function $g$ differentiable at $\theta$, with derivative $\dot g(\theta)$, the plug-in MLE $g(\hat\theta_n)$ is also asymptotically efficient, with
\[\sqrt{n}\big(g(\hat\theta_n)-g(\theta)\big)\,\xrightarrow{\mathscr{L}}\,\mathscr{N}\big(0,\,[s\dot g(\theta)]^2\big),
\] 
see, \emph{e.g.}, Lehmann and Casella (1998), p. 440. 

\medskip

In what follows, we aim at giving results on asymptotically efficient maximum likelihood estimation of linear functionals of the mixing distribution,
beyond the case of a degenerate mixing distribution localized at a point $\theta$ on the real line. As recalled in Sect.~\ref{intro}, in the deconvolution problem with the Laplace kernel density, \emph{a} NPMLE $\hat F_n$ always exists and consistency at a \emph{continuous} distribution function $F_0$ holds, but little is known about the asymptotic behaviour of the plug-in NPMLE for linear functionals. The following proposition states that, except for the above recalled degenerate case, estimation of integral linear functionals at $n^{-1/2}$-rate is impossible.
\begin{proposition}\label{prop:2}
Let $F_0$ be a non-degenerate probability measure supported on $\mathscr Y$. Let 
$\psi(F_0)$ be any integral linear functional evaluated at $F_0$.
Then, there exists no estimator sequence for $\psi(F_0)$ that is regular at $F_0$.
\end{proposition}
Some comments on Proposition \ref{prop:2}, whose proof is deferred to the Appendix B, are in order. It states that
no integral linear functional is estimable at parametric rate, in particular, by the plug-in NPMLE $\psi(\hat F_n)$. 
One can thus expect estimation, performed by any method, only at slower rates and, possibly, with a non-Gaussian limiting distribution, even if
the theorem we invoke to establish Proposition \ref{prop:2} does not give any indication about which rates to expect when 
estimation at $n^{-1/2}$-rate fails, an issue that requires further investigation. A related open question concerns the possible extension of the negative result
of Proposition \ref{prop:2} to convolution models with general kernel densities that are symmetric about zero, but not differentiable at it, a feature that
seems to play a crucial role in causing failure of estimation at parametric rate.
To sum-up, only in the case of a degenerate mixing distribution at a point $\theta$, the MLE $\hat \theta_n$ is asymptotically efficient for the location parameter and
the plug-in MLE $g(\hat\theta_n)$ is asymptotically efficient for any $g(\theta)$, with $g$ differentiable at $\theta$.




\section{Final remarks}\label{sec:3}
In this article, we have studied asymptotically efficient maximum likelihood estimation 
of linear functionals of the mixing distribution
in a standard additive measurement error 
model, when the error has either the exponential or Laplace distribution. 
In the former case, the plug-in NPMLE of certain linear functionals
is $\sqrt{n}$-consistent, asymptotically normal, efficient
and equivalent to na\"ive estimators that are empirical averages of a given 
transformation of the observations. In the latter case, instead, 
even if the kernel is generated by symmetrization about the origin of the exponential density,
left aside the degenerate case of a single Laplace model in which the MLE, the sample median, is asymptotically efficient for the location parameter, 
asymptotically efficient estimation of linear functionals completely fails, in the sense that 
estimation at $n^{-1/2}$-rate is impossible for linear functionals of non-degenerate mixing distributions.
An open question then is whether this negative result 
extends to general kernel densities symmetric about zero, but not differentiable at zero, a feature that
seems to play a crucial role in causing the failure.

\section*{Appendix A}
In this section, we present 
the proof of Proposition \ref{prop:1} on the asymptotic efficiency of the plug-in NPMLE for integral linear functionals of the mixing distribution in a convolution model with the \emph{exponential} kernel density on $[0,\,+\infty)$.

\paragraph{Proof of Proposition \ref{prop:1}}
We appeal to Theorem 2.1 of van de Geer (1997), p. 21 
(see also Theorem 11.8 of van de Geer (2000), 
pp. 217--220, for a slightly more general version) 
and, in showing that Conditions 1--4
are satisfied, we follow the indications exposed 
in Sect. 3, \emph{ibid.}, pp. 24--27.

\medskip

\noindent{\emph{Verification of Condition 1.}} (\emph{Consistency and rates}).\\[-8pt]

\noindent 
Under Assumption $(\bf{A1})$ that $M_{F_0}(1)<+\infty$, the MLE $p_{\hat F_n}$ converges in the Hellinger distance $d_{\mathrm H}$, defined as the $L^2$-distance between the square-root densities, at the rate $O_{\mathbf{P}}(n^{-1/3})$. In symbols, for $\delta_n:=n^{-1/3}$,
\[d_{\mathrm H}(p_{\hat F_n},\,p_0):=\|p_{\hat F_n}^{1/2}-p_0^{1/2}\|_2=O_{\mathbf{P}}(\delta_n).\]
The result can be obtained by applying Theorem 7.4 in van de Geer (2000), pp. 99--100, see also \emph{ibid.}, 
p. 124. 
As a consequence, see, \emph{e.g.},
Corollary 7.5, \emph{ibid.}, p. 100,
\[
\mathbb{P}_n\log \frac{2p_{\hat F_n}}{p_{\hat F_n}+p_0}
=O_{\mathbf P}(\delta_n^2),
\]
where $\delta_n^2=o(n^{-1/2})$. Consistency of $\psi(\hat F_n)$ is guaranteed by Assumption $(\bf{A2})$.


\medskip

\noindent{\emph{Verification of Condition 2.}} (\emph{Existence of the worst possible sub-directions and efficient influence functions. Differentiability of $\psi$ in a neighborhood of $F_0$}).\\[-8pt]

\noindent
For real numbers $M>M_{F_0}(1)>0$ and $r>0$, let 
\begin{equation}\label{eq:set}
\mathscr P_0:=\{F\in \mathscr P:\,\mathrm{support}\,(F)=\mathscr Y,\,\,\,M_F(1)<M,\,\,\,d_{\mathrm H}(p_F,\,p_0)\leq r\}
\end{equation}
be a Hellinger-type ball centered at $p_0$ with radius $r>0$. 
For every $\alpha\in[0,\,1)$ and $F\in\mathscr P_0$, let $F_\alpha:=\alpha F + (1-\alpha)F_0$. 
We prove
\begin{enumerate}
\item[(a)] \emph{existence of the worst possible sub-directions $h_{F_\alpha}$ such that $h_{F_\alpha}\in L^2(F_\alpha)$ and $\int h_{F_\alpha}\,\d F_\alpha=0$};\\[-8pt]
\item[(b)] \emph{existence of the efficient influence functions $b_{F_\alpha}:=A_{F_\alpha}h_{F_\alpha}$, 
where $A_{F_\alpha}h_{F_\alpha}(\cdot):=\mathrm E[h_{F_\alpha}(Y)\mid X=\cdot]$};\\[-8pt]
\item[(c)] \emph{differentiability of $\psi$ at $F_\alpha$}: 
\begin{equation}\label{eq:diff}
A^*b_{F_\alpha}(Y)=a(Y)-\psi({F_\alpha})\quad a.s.\,[F_\alpha],
\end{equation}
\emph{where $A^*b_{F_\alpha}(\cdot):=\mathrm E[b_{F_\alpha}(X)\mid Y=\cdot]$}.
\end{enumerate}

\smallskip

\noindent
For every $\alpha\in[0,\,1)$, we prove the existence of $h_{F_\alpha}$ such that the corresponding $b_{F_\alpha}=A_{F_\alpha}h_{F_\alpha}$ satisfies $A^*b_{F_\alpha}(y)=a(y)-\psi({F_\alpha})$ for $F_\alpha$-almost all $y$'s.
We proceed by first deriving the expression of $b_{F_\alpha}$ as a solution of \eqref{eq:diff} 
and then proving
the existence of the corresponding worst possible sub-direction $h_{F_\alpha}$ as required in (a) and (b). The function $b_{F_\alpha}$ has to satisfy
\begin{equation}\label{eq:infl}
A^*b_{F_\alpha}(y):=
\mathrm{E}[b_{F_\alpha}(X)\mid Y=y]=
\int_{\mathscr X}
b_{F_\alpha}(x)e^{-(x-y)}1{\{x\geq y\}}\,\d x=
a(y)-\psi(F_\alpha),
\end{equation}
where $\mathscr X=[y_{\textrm{min}},\,+\infty)$, for $F_\alpha$-almost all $y$'s.
Differentiating both sides of \eqref{eq:infl} 
with respect to $y$, we get
\begin{equation}\label{eq:der}
\int_{\mathscr X}
b_{F_\alpha}(x)e^{-(x-y)}1{\{x\geq y\}}\,\d x - b_{F_\alpha}(y)=\dot a(y).
\end{equation}
Using constraint \eqref{eq:infl} in \eqref{eq:der}, we obtain that
$a(y)-\psi(F_\alpha)-b_{F_\alpha}(y)=\dot a(y)$, whence 
$b_{F_\alpha}(y)=a(y)-\dot a(y)-\psi(F_\alpha)$. 
The solution is unique up to sets of $F_\alpha$-measure zero.
By an extension to $\mathscr X$,
$b_{F_\alpha}$ is then defined as in \eqref{eq:effinfcurv}.

\medskip

\noindent
(a) \emph{Existence of $h_{F_\alpha}\in L^2(F_\alpha)$ such that $\int h_{F_\alpha}\,\d F_\alpha=0$.}\\[5pt]
Recall that
\[p_{F_\alpha}(x):=
\int_{\mathscr Y}e^{-(x-y)}1{\{y\leq x\}}\,\d F_\alpha(y),\quad x\in\mathscr X.\]
Defined the function $I_{F_\alpha}:\,\mathscr X\rightarrow\mathbb R$ as
\[x\mapsto I_{F_\alpha}(x):=\int_{\mathscr Y}
\dot a(y)\int_{\mathscr Y}
\frac{e^{-(x-u)}}{p_{F_\alpha}(x)}1{\{u\leq y\}}
\,\d F_\alpha(u)\,1{\{y\leq x\}}\,\d y,
\]
integration by parts yields that
\[
\begin{split}
&\int_{\mathscr Y}
a(y)\frac{e^{-(x-y)}}{p_{F_\alpha}(x)}1{\{y\leq x\}}\,\d F_\alpha(y)\\
&\hspace*{2.2cm}=
\int_{\mathscr Y}
a(y)\frac{\d}{\d F_\alpha(y)}
\pt{\int_{\mathscr Y}
\frac{e^{-(x-u)}}{p_{F_\alpha}(x)}1{\{u\leq y\}}\,\d F_\alpha(u)}1{\{y\leq x\}}\,\d F_\alpha(y)
\\&\hspace*{2.2cm}=
a(y)\int_{\mathscr Y}\frac{e^{-(x-u)}}
{p_{F_\alpha}(x)}1{\{u\leq y\}}\,\d F_\alpha(u)\bigg|_{y_{\textrm{min}}}^x-
I_{F_\alpha}(x)
\\&\hspace*{2.2cm}=
a(x)\int_{\mathscr Y}\frac{e^{-(x-u)}}
{p_{F_\alpha}(x)}1{\{u\leq x\}}\,\d F_\alpha(u)
-I_{F_\alpha}(x)\\&\hspace*{2.2cm}=
a(x)-I_{F_\alpha}(x)
\end{split}
\]
because $a(y_{\textrm{min}})<+\infty$ 
and ${F_\alpha}(y_{\textrm{min}})=0$ by part (iii) of Assumption $(\bf{A3})$ combined with the fact that $F\in\mathscr P_0$. 
Analogously, since 
$\dot a_+(y_{\mathrm{min}})<+\infty$,
\[\begin{split}
&\hspace*{-0.2cm}
-\int_{\mathscr Y}
\frac{\d}{\d F_\alpha(y)}\pt{\dot a(y)
\int_{\mathscr Y} e^{-(y-u)}1{\{u\leq y\}}\,\d F_\alpha(u)}\frac{e^{-(x-y)}}{p_{F_\alpha}(x)}1{\{y\leq x\}}\,\d F_\alpha(y)\\
&\qquad=
-\dot a(y)
\int_{\mathscr Y}e^{-(y-u)}1{\{u\leq y\}}\,\d F_\alpha(u)
\times\frac{e^{-(x-y)}}{p_{F_\alpha}(x)} \bigg|_{y_{\mathrm{min}}}^x\\
&\hspace*{4.5cm}
+\int_{\mathscr Y}
\dot a(y)\int_{\mathscr Y}
\frac{e^{-(x-u)}}{p_{F_\alpha}(x)}1{\{u\leq y\}}
\,\d F_\alpha(u)\,1{\{y\leq x\}}\,\d y
\\&
\qquad=
-\dot a(x)\int_{\mathscr Y}
\frac{e^{-(x-u)}}{p_{F_\alpha}(x)}1{\{u\leq x\}}
\,\d F_\alpha(u)
+I_{F_\alpha}(x)
\\&
\qquad=
-\dot a(x)+I_{F_\alpha}(x).
\end{split}
\]
Then, defined the mapping $h_{F_\alpha}:\,\mathscr Y\rightarrow \mathbb{R}$ as
\begin{equation}\label{eq:subdir}
y\mapsto h_{F_\alpha}(y):=
\pq{a(y)-\frac{\d}{\d F_\alpha(y)}\pt{\dot a(y)\int_{\mathscr Y}e^{-(y-u)}1{\{u\leq y\}}\,\d F_\alpha(u)}-\psi({F_\alpha})
},
\end{equation}
by previous computations, we have that
\begin{eqnarray}\label{eq:ineq}
\hspace*{-0.7cm}
\forall\,x\in\mathscr X,\quad A_{F_\alpha}h_{F_\alpha}(x)&:=&\mathrm E[h_{F_\alpha}(Y)\mid X=x]\nonumber\\
&\,=&
\int_{\mathscr Y}
h_{F_\alpha}(y)
\frac{e^{-(x-y)}
}{p_{F_\alpha}(x)}1{\{y\leq x\}}\,\d F_\alpha(y)\nonumber\\
&\,=&a(x)-\dot a(x)-\psi(F_\alpha)=b_{F_\alpha}(x).
\end{eqnarray}
In order to check that $h_{F_\alpha}$ has expected value $\int h_{F_\alpha}\,\d F_\alpha=0$,
it suffices to note that, by applying twice the tower rule and 
using equalities \eqref{eq:ineq} 
and \eqref{eq:infl},
\begin{eqnarray}\label{eq:null}
\int h_{F_\alpha}\,\d F_\alpha=\mathrm E[h_{F_\alpha}(Y)]&=&
\mathrm E[\mathrm E[h_{F_\alpha}(Y)\mid X]]\nonumber\\
&=&\mathrm E[A_{F_\alpha}h_{F_\alpha}(X)]\nonumber\\&
=&\mathrm E[b_{F_\alpha}(X)]\\
&=&\mathrm E[\mathrm E[b_{F_\alpha}(X)\mid Y]]\nonumber\\
&=&\mathrm E[A^*b_{F_\alpha}(Y)]\nonumber\\
&=&\mathrm E[a(Y)-\psi(F_\alpha)]=\int_{\mathscr Y}[a(y)-\psi(F_\alpha)]\, \d F_\alpha(y)=0.\nonumber
\end{eqnarray}
\noindent
Next, we show that, for every $\alpha\in[0,\,1)$ and $F\in\mathscr P_0$,
\begin{equation*}\label{eq:sup-normdirection}
\sup_{y\in\mathscr Y}|h_{F_\alpha}(y)|<+\infty,
\end{equation*}
which implies that $h_{F_\alpha}\in L^2(F_\alpha)$. Noting that
\begin{eqnarray*}\label{eq:expresequiv}
&&\hspace*{-1cm}\frac{\d}{\d F_\alpha(y)}\pt{\dot a(y)\int_{\mathscr Y} e^{-(y-u)}1{\{u\leq y\}}\,\d F_\alpha(u)}\nonumber\\
&&\hspace*{0.5cm}=\dot a(y)e^{-y}\times\frac{\d}{\d F_\alpha(y)}\pt{\int_{\mathscr Y}e^u1{\{u\leq y\}}\,\d F_\alpha(u)}\nonumber\\
&&\hspace*{4.2cm}
+\frac{\d\big(\dot a(y)e^{-y}\big)}{\d F_\alpha(y)}\times \int_{\mathscr Y}e^u1{\{u\leq y\}}\,\d F_\alpha(u)\nonumber
\\
&&\hspace*{0.5cm}=\dot a(y)
+\frac{\d\big(\dot a(y)e^{-y}\big)}{\d F_\alpha(y)}\times \int_{\mathscr Y}e^u1{\{u\leq y\}}\,\d F_\alpha(u),
\end{eqnarray*}
we can rewrite $h_{F_\alpha}$ in \eqref{eq:subdir} as
\begin{eqnarray*}\label{eq:altern}
\hspace*{-1cm}
h_{F_\alpha}(y)
&=&
\bigg\{a(y)-\dot a(y)-\alpha
\frac{\d\big(\dot a(y)e^{-y}\big)}{\d F_0(y)}\times\frac{\d F_0(y)}{\d F_\alpha(y)}
\int_{\mathscr Y}e^u1{\{u\leq y\}}\,\d F(u)\\
&&\qquad\qquad\,\,\,\,\,\, -\,(1-\alpha)
\frac{\d\big(\dot a(y)e^{-y}\big)}{\d F_0(y)}\times\frac{\d F_0(y)}{\d F_\alpha(y)}
\int_{\mathscr Y}e^u
1{\{u\leq y\}}\,\d F_0(u)\\
&&\qquad\qquad\,\,\,\,\,\,-\,\alpha\psi(F)-(1-\alpha)\psi(F_0)\bigg\}.
\end{eqnarray*}
To conclude that $h_{F_\alpha}$ is bounded on $\mathscr Y$,
we observe two facts. First,
\[\begin{split}
|\alpha\psi(F)+(1-\alpha)\psi(F_0)|&<|\psi(F_0)|+|\psi(F)-\psi(F_0)|\\
&\leq|\psi(F_0)|+\int_{\mathscr{Y}} |a(y)|\,\d (F+F_0)(y)\\
&\leq|\psi(F_0)|+2\sup_{y\in\mathscr Y}|a(y)|<+\infty,
\end{split}
\]
where $|\psi(F_0)|<+\infty$ by Assumption 
$(\bf{A0})$ and $\sup_{y\in\mathscr Y}|a(y)|<+\infty$ by parts (i) and (ii) of Assumption 
$(\bf{A3})$. Second, for every $\alpha\in[0,\,1)$ and $F\in\mathscr P_0$,
\begin{equation}\label{eq:12}
0<\sup_{y\in\mathscr Y}\frac{\d F_0}{\d F_\alpha}(y)\leq\frac{1}{(1-\alpha)},
\end{equation}
where $(\d F_0/\d F_\alpha)$ exists because $F_\alpha$ dominates $F_0$. 
The bound in \eqref{eq:12} holds uniformly over $\mathscr P_0$. Therefore,
\begin{eqnarray}\label{eq:suph}
\hspace*{-0.8cm}
\sup_{y\in \mathscr Y}
|h_{F_\alpha}(y)|
&<&\sup_{y\in \mathscr Y}|a(y)|+
\sup_{y\in \mathscr Y}|\dot a(y)|+\frac{1}{1-\alpha}
\sup_{y\in \mathscr Y}\bigg|\frac{\d\big(\dot a(y)e^{-y}\big)}{\d F_0(y)}\bigg|[M_F(1)+M_{F_0}(1)]\\
&&\qquad\quad\,\,\,\,\,\,+\,|\psi(F_0)|+2\sup_{y\in\mathscr Y}|a(y)|<+\infty\nonumber
\end{eqnarray}
by Assumptions $(\bf{A0})$, $(\bf{A1})$, 
$(\bf{A3})$--$(\bf{A5})$ and the fact that $M_F(1)$ 
is bounded by a constant $M$ on $\mathscr P_0$.\\

\smallskip

\noindent
(b)--(c) \emph{Definition of $b_{F_\alpha}$ and differentiability of $\psi$ at $F_\alpha$.}\\[5pt]
The function $b_{F_\alpha}$ defined in \eqref{eq:effinfcurv}, which solves equation  
\eqref{eq:infl}, is such that $A_{F_\alpha}h_{F_\alpha}(x)=b_{F_\alpha}(x)$ for every $x\in\mathscr X$, in virtue of \eqref{eq:ineq}.




\medskip

\noindent
{\emph{Verification of Condition 3.}} (\emph{Control on the worst possible sub-directions $h_{F_\alpha}$}).\\[-8pt]

\noindent
Recall that $M_F(1)$ in \eqref{eq:suph}  
is bounded by $M$ on $\mathscr P_0$. Besides, the 
factor $(1-\alpha)^{-1}$, which diverges to $+\infty$ as $\alpha\rightarrow 1$, is counterbalanced by $1-\alpha$. There thus exists 
a positive constant $B\equiv B(M,\,r)<+\infty$ such that
\[\sup_{F\in\mathscr P_0}\,\sup_{0\leq\alpha<1}\,\sup_{y\in\mathscr Y}(1-\alpha)|h_{F_\alpha}(y)|\leq B.\]

\medskip

\noindent{\emph{Verification of Condition 4.}} (\emph{Control on the efficient influence functions $b_{F_\alpha}$}).\\[-8pt]

\noindent
The information for estimating $\psi(F_0)$ is positive and finite, $0<\|b_{F_0}\|_{2,P_0}^2<+\infty$. 
Also, the influence functions are uniformly bounded. In fact, for every $x\in\mathscr X$, we have 
$|b_{F_\alpha}(x)|<|b_{F_0}(x)|+|\psi(F)-\psi(F_0)|$ so that
$$\sup_{F\in\mathscr P_0}\sup_{0\leq \alpha<1}\sup_{x\in\mathscr X}|b_{F_\alpha}(x)|<3\sup_{y\in\mathscr{Y}}|a(y)|+
\sup_{y\in\mathscr{Y}}|\dot a(y)|+ |\psi(F_0)|<+\infty
$$
by Assumptions $(\bf{A0})$, $(\bf{A3})$ (parts (i) and (ii)) and $(\bf{A4})$.

\smallskip

\noindent
Next, to show that relationships (2.10) and (2.11) in van de Geer (1997), 
p. 21, are satisfied, 
we follow the reasoning illustrated in Sect. 3.4, \emph{ibid.}, pp. 26--27, 
and check that, for some positive sequence $r_n\rightarrow0$, 
\begin{equation}\label{eq:41}
\lim_{n\rightarrow+\infty}\sup_{F\in\mathscr P_n}\sup_{0\leq\alpha<1}\|b_{F_\alpha}-b_{F_0}\|_{2,P_0}^2=0,
\end{equation}
where $\mathscr P_n$ is the set obtained from $\mathscr P_0$ in \eqref{eq:set} by replacing $r$ with $r_n$. Note that $b_{F_\alpha}-b_{F_0}=\alpha[\psi(F)-\psi(F_0)]=
\alpha\int_{\mathscr{Y}}a(y)\,\d (F-F_0)(y)$. Using
integration by parts, together with conditions (i) and (ii) of Assumption $(\bf{A3})$, which jointly guarantee that $a$ is bounded on $\mathscr Y$, as well as the fact that every $F\in\mathscr P_n$ has the same support as $F_0$, we find that $\int_{\mathscr{Y}}a(y)\,\d (F-F_0)(y)=-\int_{\mathscr{Y}}\dot a(y) (F-F_0)(y)\,\d y$. 
The latter integral can be bounded above by applying 
inequality (30) in Scricciolo (2018), p. 358, which relates 
the $L^1$-Wasserstein or Kantorovich distance $W_1(F,\,F_0)=\|F-F_0\|_1$ between distribution functions $F$ and $F_0$ to the Hellinger distance between the corresponding mixtures (of exponential densities) $d_{\mathrm H}\equiv d_{\mathrm H}(p_F,\,p_0)=\|p_F^{1/2}-p_0^{1/2}\|_2$,
\begin{equation}\label{eq:wass}
W_1(F,\,F_0)\lesssim \sqrt{d_{\mathrm H}}\log^{3/4}(1/d_{\mathrm H}),
\end{equation}
where \vir{$\lesssim$} indicates inequality valid up to a constant multiple that is universal or fixed within the context, but anyway inessential for our purposes because the bound 
is uniform over $\mathscr P_n$. The inequality is obtained by setting $p=1$ and $\beta=1$, the latter value being determined by condition (29), \emph{ibid.}, p. 358, on the Fourier transform of a standard exponential density. By Assumption $(\bf{A4})$, which guarantees that $\dot a$ is bounded on $\mathscr Y$, and inequality \eqref{eq:wass}, we have
\[
\begin{split}
\|b_{F_\alpha}-b_{F_0}\|_{2,P_0}^2&<\bigg|
\int_{\mathscr{Y}}a(y)\,\d (F-F_0)(y)\bigg|^2\\
&=\bigg|
\int_{\mathscr{Y}}\dot a(y) (F-F_0)(y)\,\d y\bigg|^2\\
&\leq\bigg(\sup_{y\in\mathscr Y}|\dot a(y)|\bigg)^2W_1^2(F,\,F_0)\lesssim d_{\mathrm H}\log^{3/2}(1/d_{\mathrm H}),
\end{split}
\]
where $\lim_{n\rightarrow +\infty}d_{\mathrm H}\log^{3/2}(1/d_{\mathrm H})=0$ because $d_{\mathrm H}\leq r_n$ on $\mathscr P_n$. The limit in \eqref{eq:41} follows. 

\smallskip

\noindent
It remains to check that, for the collection of functions $\mathscr I:=\{b_{F_\alpha}:\,d_{\mathrm H}(p_F,\,p_0)\leq r, \,\,\,0\leq \alpha<1\}$, the bracketing integral
\begin{equation}\label{eq:intentropy}
\int_0^1\sqrt{\log N_{[]}\big(\varepsilon,\,\mathscr{I},\,L^2(P_0)\big)}\,\d \varepsilon<+\infty,
\end{equation}
where $N_{[]}(\varepsilon,\,\mathscr{I},\,L^2(P_0))$ is the $\varepsilon$-bracketing number of $\mathscr{I}$ for the $L^2(P_0)$-metric, namely, the smallest number of $\varepsilon$-brackets needed to cover $\mathscr{I}$, see, \emph{e.g.}, 
2.1.6 Definition (Bracketing numbers) in van der Vaart and Wellner (1996), p. 83, or Definition 2.2 in van der Geer (2000), p. 16. 
Under Assumption $(\bf{A3})$ (parts (i) and (ii)) and Assumption $(\bf{A4})$, by the same arguments as before, the $L^2(P_0)$-distance between the 
lower and upper functions $b_{F^{L}_\alpha}$ and $b_{F^{U}_\alpha}$ of every bracket $[b_{F^{L}_\alpha},\,b_{F^{U}_\alpha}]$ can be bounded above as follows:
$$\|b_{F^{U}_\alpha}-b_{F^{L}_\alpha}\|_{2,P_0}\lesssim \|F^{U}-F^{L}\|_1.$$ By 2.7.5 Theorem in
van der Vaart and Wellner (1996), pp. 159--162, the bracketing entropy of the class of all uniformly bounded, monotone functions on the real line is of the order $O(1/\varepsilon)$. Therefore, $\log N_{[]}(\varepsilon,\,\mathscr{I},\,L^2(P_0))=O(1/\varepsilon)$ and the integral in \eqref{eq:intentropy} is finite. The proof of Condition 4 is thus complete.

\medskip

\noindent
The conclusion of Theorem 2.1 follows:
\[
\psi(\hat F_n)-\psi(F_0)=\int b_{F_0}\,\d(\mathbb{P}_n-P_0) +o_{\mathbf{P}}(n^{-1/2})=
\mathbb{P}_n b_{F_0}+o_{\mathbf{P}}(n^{-1/2}),\]
where $\mathbb{P}_n b_{F_0}$ has expected value $P_0 b_{F_0}=0$, as it can be deduced from 
\eqref{eq:null} when $\alpha=0$.
Hence,
$$\sqrt{n}\big(\psi(\hat F_n)-\psi(F_0)\big)\,\xrightarrow{\mathscr{L}}\,\mathscr{N}(0,\,\|b_{F_0}\|_{2,P_0}^2)$$
and the proof is complete. \qed

\section*{Appendix B}
In this section, we present the proof of Proposition \ref{prop:2} which states that no integral linear functional of a non-degenerate
mixing distribution in a convolution model with the \emph{Laplace} kernel density
is estimable at parametric rate, in particular, by the maximum likelihood method.

\paragraph{Proof of Proposition \ref{prop:2}}
We let, at the outset, $\psi(F_0)$ be any integral linear functional,
as defined in \eqref{eq:2}, evaluated at the \vir{point} $F_0$. Arguments are laid down to identify functions $a$ (if any) whose corresponding functionals are estimable at $n^{-1/2}$-rate. To the aim, we appeal to van der Vaart's differentiability theorem, which provides a necessary and sufficient condition for pathwise differentiability of a (not necessarily linear) functional,
see Theorem 3.1, Corollaries 3.2, 3.3 and Lemma 3.4 of van der Vaart (1991), pp. 183--185, or 
Theorem 3.1, Corollaries 3.1, 3.2 and Proposition 3.1 in Groeneboom and Wellner (1992), pp. 24--28.
If differentiability of a functional fails, then, by Theorem 2.1 of van der Vaart (1991), p. 181, the functional is not estimable at $n^{-1/2}$-rate, 
see also chpt. 25 in van der Vaart (1998), pp. 358--432.
A necessary and sufficient condition for differentiability of an integral linear functional $\psi(F_0)$ is that, for $\mathscr X=\mathbb R$, there exists a function $b:\,\mathscr X\rightarrow \mathbb R$, with $b\in L^2(P_0)$, satisfying
\[\forall\,y\in\mathscr Y,\quad
\mathrm{E}[b(X)\mid Y=y]
=a(y)-\psi(F_0),
\]
explicitly,
\begin{equation}\label{eq:differentiability}
\forall\,y\in\mathscr Y,\quad 
\int_{\mathscr X}b(x)\frac{1}{2}e^{-|x-y|}\,\d x
=a(y)-\psi(F_0),
\end{equation}
where the conditional density of $X$, given $Y=y$, is $k(x-y)=e^{-|x-y|}/2$,
see \S\,7 in van der Vaart (1991), pp. 189--191, or Example 3.2 in Groeneboom and Wellner (1992), pp. 30--31. 
If an integral linear functional $\psi(F_0)$ is regularly estimable, then the condition in \eqref{eq:differentiability} must be necessarily satisfied 
and a regular estimator for $\psi(F_0)$ 
is given by $\mathbb P_nb=n^{-1}\sum_{i=1}^nb(X_i)$.
The following arguments are aimed at deriving the expression of $b$. 
Let $y\in\mathscr Y$ be fixed. For a function $a:\,\mathscr Y\rightarrow \mathbb R$
such that 
\begin{equation}\label{eq:456}
\lim_{u\rightarrow-\infty}a(u)e^u=0,
\end{equation}
where, in the case when $\mathscr Y$ is bounded, $a$ (hence its derivative $\dot a$) is taken to be identically equal to zero on $\mathscr Y^c$ so that the limit is automatically verified, integration by parts yields that
\begin{eqnarray}\label{eq:def}
\hspace*{-0.6cm}
\int_{\mathscr X}
\dot a(x)\frac{1}{2}e^{-(y-x)}1{\{x\leq y\}}\,\d x&=&
a(x)\frac{1}{2}e^{-(y-x)}\bigg|_{-\infty}^y-\int_{\mathscr X}
a(x)\frac{1}{2}e^{-(y-x)}1{\{x\leq y\}}\,\d x\\
&=&\frac{1}{2}a(y)-\int_{\mathscr X}
a(x)\frac{1}{2}e^{-(y-x)}1{\{x\leq y\}}\,\d x,\nonumber
\end{eqnarray}
whence
\begin{equation}\label{eq:1integral}
\int_{\mathscr X}[a(x)+\dot a(x)]\frac{1}{2}e^{-(y-x)}1{\{x\leq y\}}\,\d x=\frac{1}{2}a(y).
\end{equation}
The integral analogous to the one on the left-hand side of \eqref{eq:def}, but with the right branch of the Laplace density, can be dealt with similarly.
For some $a$ satisfying the limit in \eqref{eq:456} and also
\begin{equation*}\label{eq:378}
\lim_{u\rightarrow +\infty}a(u)e^{-u}=0,
\end{equation*}
where the same proviso on $a$ and $\dot a$ applies for the case when
$\mathscr Y$ is bounded, we get 
\[\begin{split}
\int_{\mathscr X}\dot a(x)\frac{1}{2}e^{-(x-y)}1{\{x>y\}}\,\d x&=
a(x)\frac{1}{2}e^{-(x-y)}\bigg|_{y}^{+\infty}+
\int_{\mathscr X}a(x)\frac{1}{2}e^{-(x-y)}1{\{x>y\}}\,\d x\\
&=-\frac{1}{2}a(y)+\int_{\mathscr X}
a(x)\frac{1}{2}e^{-(x-y)}1{\{x>y\}}\,\d x,
\end{split}
\]
whence
\begin{equation}\label{eq:2integral}
\int_{\mathscr X}[a(x)-\dot a(x)]\frac{1}{2}e^{-(x-y)}1{\{x>y\}}\,\d x=\frac{1}{2}a(y).
\end{equation}
Summing side by side \eqref{eq:1integral}
and \eqref{eq:2integral} and subtracting $\psi(F_0)$ on both sides of the resulting equation, we obtain
\[
\int_{\mathscr X}
[a(x)-\mathrm{sgn}(x-y)\dot a(x)-\psi(F_0)]\frac{1}{2}e^{-|x-y|}\,\d x=a(y)-\psi(F_0).
\]
In order to get rid of the dependence of the function $a(\cdot)-\mathrm{sgn}(\cdot-y)\dot a(\cdot)-\psi(F_0)$ on $y$, the derivative $\dot a(\cdot)$ must be equal to zero, which means that $a(\cdot)$ is identically equal to a constant on ${\mathscr Y}$ and the functional is trivially equal to the constant itself. 
Conclude that there exists no integral linear functional $\psi(F_0)$ of a non-degenerate mixing distribution $F_0$ that can be estimated at $n^{-1/2}$-rate. This completes the proof.
\qed

\section*{References}
\begin{small}
\begin{description}
%
%


\item[]
Billingsley P (1995) Probability and measure.
Wiley, New York, 3rd edition

\item[]
Bolthausen E, Perkins E, van der Vaart A (2002)
Lectures on probability theory and statistics. Ecole d'Et\'e de Probabilit\'es
de Saint-Flour XXIX -- 1999. Bernard P (ed) Lecture Notes in Mathematics, Vol 1781.
Springer-Verlag, Berlin, pp 331--457



\item[]
Buonaccorsi, JP (2010) Measurement error: models, methods, and applications. Chapman \& Hall/CRC Press, Boca Raton, FL

\item[]
Buzas JS, Stefanski LA, Tosteson TD (2005) Measurement error. 
In: Ahrens W, Pigeot I (eds) Handbook of epidemiology. 
Springer-Verlag, Berlin, Heidelberg, pp 729--765

\item[]
Carroll LB (2017) Nuclear steam generator fitness-for-service assessment. In: Riznic J (ed)
Steam generators for nuclear power plants.
Woodhead Publishing, pp 511--523

\item[]
Carroll RJ, Hall P (1988) Optimal rates of convergence for
deconvolving a density. J Amer Statist Assoc 83:1184--1186 


\item[]
Chen J (2017) Consistency of the MLE under mixture models. Stat Sci 32:47--63


\item[]
Daniels HE (1961) The asymptotic efficiency of a maximum likelihood estimator. In: Proc Fourth Berkeley Symp on Math Statist and Prob, Vol 1. Univ of Calif Press, pp 151--163



\item[]
Dattner I, Goldenshluger A, Juditsky A (2011) On deconvolution of distribution functions. 
Ann Stat 39:2477--2501

\item[]
Davidian M, Lin X, Morris JS, Stefanski LA (2014)
The work of Raymond J. Carroll: The impact and influence of a statistician.
Springer International Publishing, Switzerland



\item[]
Easterling RG (1980) Statistical analysis of steam generator inspection plans and eddy current testing. Washington, DC: Division of Operating Reactors, Office of Nuclear Reactor Regulation, US Nuclear Regulatory Commission

\item[]
Fan J (1991) On the optimal rates of convergence for nonparametric
deconvolution problems. Ann Stat 19:1257--1272

\item[]
Fuller WA (1987) Measurement error models. 
John Wiley, New York

\item[]
Groeneboom P, Wellner JA (1992) Information bounds and nonparametric maximum likelihood estimation. Birkh\"auser, Basel


\item[]
Hall P, Lahiri SN (2008) Estimation of distributions, moments and quantiles in deconvolution problems. 
Ann Stat 36:2110--2134

\item[]
Huber PJ (1967) The behavior of maximum likelihood estimates under nonstandard conditions. In: Proc Fifth Berkeley Symp on Math Statist and Prob, Vol 1.
Univ of Calif Press, pp 221--233


\item[]
Kotz S, Kozubowski TJ, Podg\'orski K (2001) The Laplace distribution and generalizations:
a revisit with applications to communications,
economics, engineering, and finance. 
Birkh\"auser, Boston



\item[]
Laplace P-S (1774) M\'emoire sur la probabilit\'e des causes par les \'ev\'enements. 
M\'em Acad Roy Sci Paris (Savants \'etrangers) 
Tome VI:621--656

\item[]
Lehmann EL, Casella G (1998) Theory of point estimation, 2nd ed. Springer-Verlag, New York

\item[]
Lindsay BG (1983) The geometry of mixture likelihoods: a general theory. Ann Stat 11:86--94



\item[]
Norton RM (1984) The double exponential distribution: using calculus to find a maximum likelihood estimator. Am Stat 38:135--136


\item[]
Scricciolo C (2018) Bayes and maximum likelihood for $L^1$-Wasserstein deconvolution of Laplace mixtures. Stat Methods Appl 27:333--362

\item[]
Sollier T (2017) Nuclear steam generator inspection and testing. In: Riznic J (ed)
Steam generators for nuclear power plants.
Woodhead Publishing, pp 471--493

\item[]
Stefanski L, Carroll RJ (1990) Deconvoluting kernel density
estimators. Statistics 21:169--184 



\item[]
van de Geer S (1997) Asymptotic normality in mixture models. ESAIM Probab Stat 1:17--33

\item[]
van de Geer SA (2000) Empirical processes in M-estimation. Cambridge University Press, Cambridge

\item[]
van de Geer S (2003) Asymptotic theory for maximum likelihood in nonparametric mixture models. Comput Stat Data An
 41:453--464

\item[]
van der Vaart A (1991) On differentiable functionals.
Ann Stat 19:178--204


\item[]
van der Vaart AW (1998) Asymptotic statistics. 
Cambridge University Press, Cambridge


\item[]
van der Vaart AW, Wellner JA (1996) Weak convergence and empirical processes. Springer-Verlag, New York

\item[]
Vardi Y (1989) Multiplicative censoring, renewal processes, deconvolution
and decreasing density: nonparametric estimation. Biometrika 76:751--761

\end{description}
\end{small}










\end{document}